\documentclass[11pt,a4paper,reqno]{amsart}

\usepackage{fullpage}
\usepackage{amssymb}
\usepackage[dvipsnames]{xcolor}
\usepackage{hyperref}
\hypersetup{
pdfauthor = {Hughes+Israel+Mayeli}, 
pdftitle = \text{Trace bounds on limiting operators for rough domains}, 
colorlinks = true, 
linkcolor = blue, 
citecolor = red
}

\numberwithin{equation}{section}

\theoremstyle{plain}
\newtheorem{theorem}{Theorem}[section]
\newtheorem{lemma}[theorem]{Lemma}
\newtheorem{corollary}[theorem]{Corollary}
\newtheorem{proposition}[theorem]{Proposition}

\theoremstyle{definition}
\newtheorem{definition}[theorem]{Definition}

\theoremstyle{remark}
\newtheorem{remark}[theorem]{Remark}

\newcommand{\Z}{\mathbb Z}
\newcommand{\N}{\mathbb N}
\newcommand{\R}{\mathbb{R}}

\newcommand{\cM}{\mathcal M}

\newcommand{\dist}{\mathrm{dist}}
\newcommand{\diam}{\mathrm{diam}}
\newcommand{\tr}{\mathrm{tr}}

\DeclareMathOperator{\Per}{Per}
\DeclareMathOperator{\HS}{HS}
\newcommand{\1}{\mathbf{1}}

\def\dd{{\,{\rm d}}}

\begin{document}
\title{Trace bounds for  limiting operators on rough domains}

\author{Kevin Hughes}

\author{Arie Israel}
\thanks{A.~ Israel was supported by the Air Force Office of Scientific Research, under award FA9550-19-1-0005 and the National Science Foundation grant DMS-2453770.}

\author{Azita Mayeli}
\thanks{A.~Mayeli was supported in part by the National Science Foundation grant DMS-2453769, the AMS-Simons Research Enhancement Grant, and the PSC-CUNY grants 67807-00 56 and 67807-00 57.}
\date{\today}


\begin{abstract} 
This work concerns a quantitative form of Landau's eigenvalue theorem for spatio-spectral limiting operators.  
We isolate a simple mechanism that converts the problem of estimating the distribution of eigenvalues of a limiting operator into the problem of bounding the trace of the difference between the operator and its square. 
This mechanism allows us to analyze limiting operators for domains with fractal boundaries. 
When the boundaries have finite perimeter, we recover the expected optimal dependence on the scaling parameter.
\end{abstract}

\maketitle


\section{Introduction}

Let \(d \in \N\) denote the dimension of the ambient Euclidean space. 
Let \(\mathcal{F}\) denote the Fourier transform on \(L^2(\R^d)\) and \(\mathcal{F}^{-1}\) denote its inverse; both are \(L^2(\R^d)\) normalized so that 
\begin{equation}\label{eq:Fourier-convention}
\mathcal{F}f(\xi) := (2\pi)^{-d/2} \int_{\R^d} f(x)e^{ix \cdot \xi} \dd{x}, 
\quad 
\mathcal{F}^{-1}g(x) := (2\pi)^{-d/2} \int_{\R^d} g(\xi)e^{-ix \cdot \xi} \dd{\xi} 
.\end{equation}
Let \( F,S \subset \R^d\) be Lebesgue measurable sets. 
Define the orthogonal spatial projection \(P_F\) on \(L^2(\R^d)\) by \(P_F f := \mathbf{1}_F \cdot f\). Define the band-limiting projection \(B_S\) on \(L^2(\R^d)\) by \(B_S := \mathcal{F}^{-1}P_S\mathcal{F}\). Compose these projections to form the spatio-spectral limiting operator \(T_{F,S} = P_F B_S P_F\) acting on \(L^2(\R^d)\). This operator acts to limit a function $f(x)$ to the region $F$ in the physical domain, and to the region $S$ in the Fourier domain.\footnote{We could have instead defined the operator as the product \(B_S P_F B_S\), but both expressions have the same non-zero spectrum.} If both \(F\) and \(S\) have finite Lebesgue measure, \(T_{F,S}\) is compact and self--adjoint on $L^2(\R^d)$, with eigenvalues ordered so that 
\(
1 \geq \lambda_1(T_{F,S}) \geq \lambda_2(T_{F,S}) \geq \cdots \geq 0.
\)

Given a subset \(U\) in \(\R^d\) and \(R>0\), define the isotropic dilation \(RU:= \{Rx : x \in U\}.\) 
For any Lebesgue measurable set \(U\) in \(\R^d\), let \(|U|\) denote its Lebesgue measure. 
When \(F,S\) are fixed, define the spatio-spectral limiting operators \(T_R := T_{RF,S} = P_{RF} B_S P_{RF}\) for each \(R>0\). 
For \(\varepsilon \in (0,1)\), let \(N_\varepsilon(R) := \#\{k:\lambda_k(T_R)>\varepsilon\}\) denote the number of eigenvalues of \(T_R\) that are greater than \(\varepsilon\). 
A classical result of Landau \cite{Landau75} shows that
\begin{equation}\label{Landau:asymptotic}
\lim_{R \to \infty} R^{-d}N_\varepsilon(R) = (2\pi)^{-d}|F|\cdot|S|.
\end{equation}
The right-hand side of \eqref{Landau:asymptotic} gives the leading-order number of effective degrees of freedom for signals localized in space \(F\) and frequency \(S\). This asymptotically confirmed the classical engineering use of the Nyquist rate and the \(2WT\) Theorem (also known as the Sampling Theorem).

Since \eqref{Landau:asymptotic} only requires that the sets \(F,S\) have finite and positive Lebesgue measure, the original argument in~\cite{Landau75} does not provide quantitative rates. In dimension one, Landau and Widom~\cite{LandauWidom80} proved a sharper two-term asymptotic describing the logarithmic transition region for time- and frequency-limiting operators when both sets are a finite union of intervals. In higher dimensions, quantitative estimates are sensitive to the geometry of the boundaries of \(F\) and \(S\) and have been studied under various regularity assumptions on the underlying domains. See \cite{israelmayeli2023acha,MaRoSp23,hughes2024eigenvalue,HIM-disk} and the recent development \cite{Kulikov-Larsen}.

Our goal is to present a simple quantitative variant of Landau's proof \cite[Theorem~1]{Landau75} that naturally incorporates quantitative geometric information about the boundaries of \(F\) and \(S\). Our ensuing setup allows us to simultaneously quantify the plunge region by choosing a different function \(\varphi\) in the following lemma. As in \cite{Landau75}, the main device for doing so is the following reduction to bounding the trace of \(T_{F,S}-T_{F,S}^2\). 

\begin{lemma}\label{lemma:trace_bound}
Let \(T\) be a positive trace-class contraction. 
If \(\varphi:[0,1]\to \mathbb C\) is a measurable function satisfying \(\varphi(0)=0\) and
\[
C_\varphi
:=
\sup_{0<t<1} \frac{|\varphi(t)-\varphi(1)t|}{t(1-t)}
<
\infty,
\]
then $\varphi(T)$ is trace-class and
\[
|\operatorname{tr}(\varphi(T))-\varphi(1)\operatorname{tr}(T)|
\leq
C_\varphi \tr(T-T^2).
\]
\end{lemma}

We use this by analyzing the following identity for the trace of \(T_{F,S}-T_{F,S}^2\). 
\begin{proposition}\label{prop:A-fxn}
If $F, S\subset \R^d$ are sets of finite Lebesgue measure, then 
\begin{equation}\label{id:D-fxn}
\tr(T_{F,S}-T_{F,S}^2)
= (2\pi)^{-d} \int_{\R^d}|\mathcal{F}^{-1}(\1_S)(z)|^2 \cdot |F\cap(F^c-z)| \dd{z}
.\end{equation}
\end{proposition}

Our main result bounds \(\tr(T_{F,S}-T_{F,S}^2)\) when our sets \(F,S\) satisfy the following condition. 
\begin{definition}\label{def:perimeter}
Let \(0<\gamma \leq 1\). 
For \(U\subseteq \R^d\) a Lebesgue measurable set of finite and positive measure, define the \emph{$\gamma$-perimeter} of $U$ by
\[
\Per_\gamma(U)
:= \sup_{h \in \R^d \setminus \{0\}} \frac{|U\triangle(U-h)|}{|h|^\gamma}.
\]
When \(U\) has finite \(\gamma\)-perimeter, define the scale \( \delta_U :=  \left( \frac{|U|}{\Per_\gamma(U)} \right)^{1/\gamma}\).
\end{definition}

For \(0<\gamma \leq 1\), the finiteness of \(\Per_\gamma(U)\) is equivalent to \(\1_U\in \dot{B}^\gamma_{1,\infty}(\R^d)\),  the homogeneous Besov space.
At the endpoint \(\gamma=1\), the condition \(\Per_\gamma(U)<\infty\) is equivalent to \(\mathbf 1_U\in \dot{BV}(\R^d)\), the space of functions of bounded variation. 
Thus, \(\Per_1(U)<\infty\) is equivalent to \(U\) having finite perimeter in the usual sense. See \cite{sickel2020regularity, triebel1983theory} for more information on these equivalences, related discussions, and examples. 

\begin{theorem}\label{thm:trace-of-defect-bd:Besov}
Let \(0<\gamma,\eta\leq 1\). 
If \(F, S\subseteq \R^d\) are sets of finite Lebesgue measure such that \(F\) has finite \(\gamma\)-perimeter and \(S\) has finite \(\eta\)-perimeter, then 
\begin{equation}\label{bd:D:symdiff}
\tr(T_{F,S} - T_{F,S}^2) 
\leq
C \cdot 
\begin{cases}
\Per_\gamma(F) \Per_\eta(S) \log_{+} (\delta_F \delta_S) &\text{if }\gamma=\eta \\
\Per_\eta(S) \Per_\gamma(F)^{\eta/\gamma} |F|^{1-\eta/\gamma} &\text{if }\gamma > \eta \\
\Per_\eta(S)^{\gamma/\eta} \Per_\gamma(F) |S|^{1-\gamma/\eta} &\text{if }\gamma < \eta
.\end{cases}
\end{equation}
Here, $C$ is a constant depending on $d, \eta$, and $\gamma$, while $\log_{+}(a) := \max\{1,\log(a)\}$.
\end{theorem}
\noindent All logarithms in this paper are base 2 unless specified otherwise.

\begin{remark}
Suppose $F$ is a bounded measurable set. Take $h \in \R^d$ with $|h| > \diam(F)$ and observe that $F$ and $F - h$ do not intersect. Thus, $| F \triangle (F-h)| = 2|F|$. This implies $\Per_\gamma(F) \geq 2|F| \diam(F)^{-\gamma}$, which yields $\delta_F \leq 2^{1/\gamma} \diam(F)$. Similarly, $\delta_S \leq 2^{1/\eta} \diam(S)$ when \(S\) is a bounded measurable set.
\end{remark}

If we scale the domain $F$ to $RF$ for $R \geq 2$, then, for the scaled family of operators $T_R = T_{RF,S}$, the estimate \eqref{bd:D:symdiff} can be expressed as 
\begin{align}\label{trace-defect-estimate}
\tr(T_R - T_R^2) 
\leq
C(F,S)
\begin{cases}
R^{d-\gamma} \log(R) &\text{if }\gamma=\eta \\
R^{d-\eta} &\text{if }\gamma > \eta \\
R^{d-\gamma} &\text{if } \gamma <\eta.
\end{cases}
\end{align}
Here, $C(F,S)$ is a constant determined by $d,\gamma,\eta$ and the perimeters and measures of $F,S$.

Observe that \(\varphi := \1_{(\varepsilon,1]}\) satisfies the hypotheses of Lemma~\ref{lemma:trace_bound} with \(C_{\varphi} \leq (\varepsilon(1-\varepsilon))^{-1}\), \(\tr(\1_{(\varepsilon,1]}(T_R)) = N_{\varepsilon}(R)\), and by Mercer's theorem, 
\begin{equation}\label{id:trace}
\tr(T_R)
= (2\pi)^{-d} |RF|\cdot|S| = (2\pi)^{-d} R^d|F|\cdot|S|
.\end{equation}
Plugging \eqref{trace-defect-estimate} and \eqref{id:trace} into Lemma~\ref{lemma:trace_bound}, we obtain the following quantitative form of \eqref{Landau:asymptotic}.

\begin{corollary}
\label{cor:quantitative-Landau}
Let \(F, S\subseteq \R^d\) be sets of finite Lebesgue measure. If \(F\) has finite \(\gamma\)-perimeter and \(S\) has finite \(\eta\)-perimeter for \(0<\gamma,\eta\leq 1\), then for all \(R\geq 2\) and \(\varepsilon\in(0,1)\),
\[
\big| N_\varepsilon(R) - (2\pi)^{-d}|F||S|R^d \big|
\leq
\frac{C(F,S)}{\varepsilon(1-\varepsilon)}
\begin{cases}
R^{d-\gamma} \log(R) &\text{if }\gamma=\eta \\
R^{d-\eta} &\text{if }\gamma > \eta \\
R^{d-\gamma} &\text{if }\eta > \gamma
.\end{cases}
\]
\end{corollary}

\begin{remark}
Define the \emph{large-scale $(d-\gamma)$-dimensional upper Minkowski content} of $E \subset \R^d$ by
\begin{equation}\label{def:UMC}
\cM^{d-\gamma}(E) := \sup_{0< r \leq \diam(E)} \left| E_r \right| r^{-\gamma},
\end{equation}
where $E_r := \{ x \in \R^d : \dist(x,E) \leq r\}$ is the $r$-neighborhood of $E$. Finiteness of $\cM^{d-\gamma}(E)$ is equivalent to finiteness of the usual $(d-\gamma)$-dimensional upper Minkowski content of $E$. 

If $U\subset \mathbb{R}^d$ is a bounded measurable set, then there exists a dimensional constant \(C_d\) such that
\begin{equation}\label{eqn:per_mink_rel}
\Per_\gamma(U) \leq C_d \cM^{d-\gamma}(\partial U).
\end{equation}
To see this, let $\Delta := \diam(U) = \diam(\partial U)$. If $|w| \leq \Delta$, then $U \triangle (U - w) \subset (\partial U)_{|w|}$ and 
\[
|U \triangle (U-w)| \cdot |w|^{-\gamma} \leq | (\partial U)_{|w|}| \cdot |w|^{-\gamma} \leq \cM^{d-\gamma}(\partial U).
\]
If $|w| > \Delta$, then 
\[
|U \triangle (U - w)| \cdot |w|^{-\gamma}
\leq 2 |U| \Delta^{-\gamma} \lesssim_d \Delta^{d-\gamma}
\lesssim_d |(\partial U)_\Delta| \Delta^{-\gamma}
\leq \cM^{d-\gamma}(\partial U),
\]
where the last inequality follows because $r=\Delta$ is an admissible scale in \eqref{def:UMC}. Combining these two estimates yields \eqref{eqn:per_mink_rel}. 

In view of \eqref{eqn:per_mink_rel}, finite $(d-\gamma)$-dimensional upper Minkowski content of the boundary implies finite $\gamma$-perimeter for bounded measurable sets. In particular, Corollary \ref{cor:quantitative-Landau} yields new eigenvalue estimates for limiting operators on domains with fractal boundaries such as the interior of the Koch snowflake in $\R^2$, which has finite $(\log_3{4})$-dimensional upper Minkowski content, and hence has finite (\(2-\log_{3}{4}\))-perimeter. 
\end{remark}

The estimate on $N_\varepsilon(R)$ implies an estimate on the size of the plunge region of $T_R$.\footnote{Alternatively derive \eqref{eqn:plunge1} from observing \(\varphi:=\1_{(\varepsilon,1-\varepsilon)}\) satisfies the hypotheses of Lemma~\ref{lemma:trace_bound} and applying \eqref{trace-defect-estimate} as in the derivation of Corollary~\ref{cor:quantitative-Landau}.} 
Let \(M_\varepsilon(R) := \#\{\lambda_k(T_R)\in(\varepsilon,1-\varepsilon)\}\) count the number of eigenvalues of $T_R$ lying in $(\varepsilon,1-\varepsilon)$ for $\varepsilon\in(0,\frac{1}{2})$. 
If $F$ and $S$ have finite $\gamma$-perimeter, then for each fixed $\varepsilon \in (0,1/2)$,
\begin{equation}\label{eqn:plunge1}
M_\varepsilon(R)
\leq N_{\varepsilon}(R)-N_{1-\varepsilon}(R)
\leq C \varepsilon^{-1} R^{d-\gamma}\log (R).
\end{equation}
Recent works such as \cite{MaRoSp23,israel15eigenvalue,Karnik21,israelmayeli2023acha,hughes2024eigenvalue,HIM-disk,Kulikov-Larsen} improve the dependence on \(\varepsilon\) under various geometric conditions on \(F\) and \(S\). 
When \(F,S\) are bounded sets whose boundaries have finite $(d-1)$-dimensional upper Minkowski content, Theorem~1.6 in \cite{Kulikov-Larsen} gives, for fixed \(\varepsilon\in(0,1/2)\), the bound 
\begin{equation}\label{eqn:KL_plunge}
M_\varepsilon(R)\leq C_\varepsilon R^{d-1}\log(R)^2
\end{equation}
where the constant \(C_\varepsilon\) scales like a power of \(\log(1/\varepsilon)\).  In the fixed-\(\varepsilon\) regime, our bound \eqref{eqn:plunge1} with $\gamma=1$ has one fewer power of \(\log (R)\), at the expense of a less sharp dependence on \(\varepsilon\). Further, \eqref{eqn:plunge1} requires only the finiteness of the perimeter of $F$ and $S$, whereas \eqref{eqn:KL_plunge} requires the stronger Minkowski condition on the boundaries of $F$ and $S$.

\section{Proof of Lemma~\ref{lemma:trace_bound}}

Let \(1 \geq  \lambda_1(T) \geq \lambda_2(T) \geq \dots > 0\) denote the positive eigenvalues of \(T\). Since \(T\) is trace-class, \(\sum_k \lambda_k(T)=\tr(T)<\infty\). 
The assumption on \(\varphi\) implies
\[
|\varphi(\lambda_k(T))-\varphi(1)\lambda_k(T)|
\leq
C_\varphi \lambda_k(T)(1-\lambda_k(T)).
\]
To complete the proof,  sum over \(k\) to obtain
\[
\left| \tr(\varphi(T)) - \varphi(1)\tr(T) \right|
\leq
C_\varphi \sum_k \lambda_k(T)(1-\lambda_k(T))
=
C_\varphi \tr(T-T^2).
\]
\section{Proof of Proposition~\ref{prop:A-fxn}}
Observe that
\(
T_{F,S}-T_{F,S}^2 = P_{F} B_S(I-P_{F})B_S P_{F}.
\)
Since $F$ and $S$ have finite measure, $B_S P_{F}$ is Hilbert--Schmidt and so is $(I-P_{F})B_S P_{F}$. 
By the cyclic property of the trace,
\[
\tr(T_{F,S}-T_{F,S}^2)
= \mathrm{tr}(P_{F} B_S(I-P_{F})B_S P_{F})
= \|(I-P_{F})B_S P_{F}\|_{\HS}^2
\]
where \(\|\cdot\|_{\HS}\) denotes the Hilbert--Schmidt norm. 
The operator $(I-P_{F})B_S P_{F}$ has kernel
\[
\mathbf{1}_{F^c}(x) (2\pi)^{-d/2}\mathcal{F}^{-1}(\1_S)(x-y) \mathbf{1}_{F}(y).
\]
Therefore,
\[
\tr(T_{F,S}-T_{F,S}^2)
= (2\pi)^{-d} \int_{y\in F} \int_{x\in F^c} \hspace{-5mm} |\mathcal{F}^{-1}(\1_S)(x-y)|^2 \dd{x} \dd{y}
.\]
Make the change of variables \(x \mapsto z:=x-y\), apply Fubini's theorem, and integrate with respect to \(y\) to complete the proof.

\section{Proof of Theorem~\ref{thm:trace-of-defect-bd:Besov}}

We first consider the case $0<\delta_F\delta_S<2$. 
If $\gamma=\eta$, then by the definition of $\delta_F$ and $\delta_S$,
\[
|F||S|
=
\Per_\gamma(F)\Per_\gamma(S)(\delta_F\delta_S)^\gamma
\leq
2^\gamma \Per_\gamma(F)\Per_\gamma(S).
\]
Using the trivial bounds $\tr(T_{F,S}-T_{F,S}^2) \leq \tr(T_{F,S}) = (2\pi)^{-d} |F| |S|$ and $\log_+(\delta_F\delta_S)\geq 1$, this gives
\[
\tr(T_{F,S}-T_{F,S}^2)
\lesssim
\Per_\gamma(F)\Per_\gamma(S) \log_+(\delta_F\delta_S).
\]
If $\gamma>\eta$, then $ 
|F||S| \lesssim \Per_\eta(S)\Per_\gamma(F)^{\eta/\gamma}|F|^{1-\eta/\gamma}$. Similarly, if $\gamma<\eta$, then $|F||S|
\lesssim \Per_\eta(S)^{\gamma/\eta}\Per_\gamma(F)|S|^{1-\gamma/\eta}$. 
This completes the proof of the estimate when $\delta_F\delta_S<2$.

Next, we treat the case $\delta_F\delta_S\geq 2$. We start with the following observation. 
Since a characteristic function only takes the values $0$ or $1$, \(\Per_\eta(S) < \infty \iff \1_S \in \dot{B}^{\eta}_{1,\infty} \iff \1_S \in \dot{B}^{\eta/2}_{2,\infty}\). Further, $\Per_\eta(S)$ is equivalent to $\| \1_S\|^2_{\dot{B}^{\eta/2}_{2,\infty}}$. 
Using the Fourier characterization of Besov spaces as in \cite[Sections 2.5.7 and 2.5.12]{triebel1983theory}, there exists a dimensional constant \(C_d\) such that 
\begin{equation}\label{bd:Fourierdecay}
\int_{2^j \leq |w| \leq 2^{j+1}} |\mathcal{F}^{-1}(\1_S)(w)|^2 \dd{w}
\leq C_d \Per_\eta(S) 2^{-j\eta} 
\quad \text{for all} \quad j \in \Z
.\end{equation}
Given a scale $\delta > 0$, sum estimate \eqref{bd:Fourierdecay} over $j \geq \log(\delta)$ to obtain
\begin{equation}\label{bd:Fourierdecay2}
\int_{|w| \geq \delta} |\mathcal{F}^{-1}(\1_S)(w)|^2 \dd{w}
\lesssim_{d,\eta} \Per_\eta(S) \delta^{-\eta} 
\quad \text{for all} \quad \delta > 0.
\end{equation}
If $j \ll 0$, the bound \eqref{bd:Fourierdecay} is weaker than that derived by Plancherel,
\(
\| \mathcal{F}^{-1}\1_S \|_{L^2}^2 = \| \1_S \|_{L^2}^2 = |S|.
\)

The perimeter condition on $F$ implies
\[
|F \cap (F^c-z)| \leq |F\triangle(F-z)| \leq \Per_\gamma(F) |z|^\gamma.
\]
Alternatively, by the volume bound on $F$,
\(
|F \cap (F^c-z)| \leq |F|.
\) 
Define scales $\delta_F$ and $\delta_S$ as in the statement of Theorem~\ref{thm:trace-of-defect-bd:Besov}. 
The first estimate is sharper when $|z| < \delta_F$, and the second estimate is sharper when $|z| > \delta_F$. Recall \eqref{id:D-fxn} from Proposition~\ref{prop:A-fxn}, split the integral into dyadic pieces, and use the above estimates to find
\begin{align*}
\tr(T_{F,S}-T_{F,S}^2)
&\leq 
\int_{|z| < 2 \delta_S^{-1}} |\mathcal{F}^{-1}(\1_S)(z)|^2 \Per_\gamma(F)|z|^\gamma  \dd{z} \\
& \qquad + 
\sum_{\log(\delta_S^{-1}) \leq j \leq \log(\delta_F)} \int_{2^j \leq |z| \leq 2^{j+1}} |\mathcal{F}^{-1}(\1_S)(z)|^2 \Per_\gamma(F)|z|^\gamma \dd{z}
\\& \qquad + \int_{|z| > \delta_F}  |\mathcal{F}^{-1}(\1_S)(z)|^2 |F| \dd{z}.
\end{align*} 
For the first integral over $|z| < 2\delta_S^{-1}$ we use Plancherel's theorem, and for the remaining integrals we use \eqref{bd:Fourierdecay} and \eqref{bd:Fourierdecay2}, to find that 
\begin{align*}
\tr(T_{F,S}-T_{F,S}^2)
&\lesssim_{d,\eta} 
|S| \Per_\gamma(F)\delta_S^{-\gamma}
+ \hspace{-3mm} \sum_{\log(\delta_S^{-1}) \leq j \leq \log(\delta_F)} \hspace{-8mm} \Per_\eta(S) \Per_\gamma(F) 2^{j(\gamma- \eta)} 
+ |F| \Per_\eta(S) \delta_F^{-\eta}.
\end{align*}
The middle sum is non-empty because $\delta_F\delta_S\geq 2$. Summing the geometric series and substituting the definitions of $\delta_F$ and $\delta_S$ completes the proof.





\end{document}